\documentclass[11pt,
]{article}

\usepackage{amssymb,amsfonts,amstext,amsmath,amsthm,
latexsym,mathrsfs,amsbsy}

\usepackage[bbgreekl]{mathbbol}

\usepackage{marginnote}
\usepackage{makeidx,mparhack}

\usepackage{graphicx}
\usepackage{xcolor}
\usepackage{doi}

\input{pa.sty}
\makeatletter
\if@twoside%
\else\fi%
\makeatother  
\begin{document}

\title{Jensen $\id13$ reals by means of ZFC$^-$ or
second order Peano arithmetic}

\author 
{Vladimir Kanovei \thanks
{Institute for Information Transmission Problems (IITP RAS), 
{\tt kanovei@googlemail.com}. 
%Partial support of RFBR Grant 17-01-00705 acknowledged.
}
}
\date{\today}
\maketitle

%\renek{\abstractname}     {Abstract}

\begin{abstract}
%Suppose that $n\ge3$. 
It was established by Jensen in 1970 
%(and also by 
%Jensen and Solovay by a somewhat different method), 
that there is a generic extension $\rL[a]$ of the 
constructible universe $\rL$ by a real $a\nin\rL$ 
such that $a$ is $\id13$ in $\rL[a]$.
Jensen's forcing construction has found a number of
applications in modern set theory.
A problem has been recently discussed whether
Jensen's construction can be reproduced entirely within
second order Peano arithmetic or equivalently ZFC$^-$
(minus the Power Set axiom).
The obstacle is that the proof of the key CCC property
(whether by Jensen's original argument or a later
proof using $\Diamond$) essentially involve 
countable elementary submodels of $\rL_{\om_2}$,
which is way beyond ZFC$^-$.
We show how to circumwent this difficulty
by means of killing only definable
antichains in the course of a Jensen-like transfinite
construction of the forcing, and then define a
model with a minimal $\varPi^1_2$ singleton as a
class-forcing extension of a model of ZFC$^-$ plus
$\rV=\rL$.
\end{abstract}

\cha{Introduction}

The next theorem is the main result of this paper.

\bte
%[$\rV=\rL$]
\lam{gt}
There exists a forcing notion\/ $J\sq\Lomi$
satisfying the following$:$
\ben
\renu
\itlb{gt1}
$J$ is parameter-free\/ $\id\Lomi1\,;$
and CCC in\/ $\Lomi$, in which it is true that$:$

\itlb{gt2}
$J$ is CCC \poo\ all antichains\/ $A\sq\Lomi$,
definable in\/ $\Lomi$ with parameters$;$

\itlb{gt3}
$J$ as a class-forcing adjoins a real\/ $a\in\dn$
to\/ $\Lomi$, such that$:$

\itlb{gt4} 
$\ans a$ is a nonconstructible\/ $\ip12$ singleton
in\/ $\Lomi[a];$ 

\itlb{gt5} 
$a$ is minimal over\/ $\Lomi$ in the sense that\/ 
$a\nin\Lomi$ but any real\/ $x\in\Lomi[a]$
either belongs to\/ $\Lomi$ or satisfies\/ 
$a\in\Lomi[x]\;.$  
\een
In addition, {\rm(*)} 
claims \ref{gt1}--\ref{gt5} are provable in\/
$\zfcm$,
and hence the theory\/ $\zfcm$ plus ``there exists
a nonconstructible minimal\/ $\ip12$ real singleton''
is equiconsistent with $\zfcm$ and thereby with
second-order Peano arithmetic $\pad$.
\ete

We recall that $\zfcm$ is the theory
\index{zzZFCm@$\zfcm$}%
\index{theory!ZFCm@$\zfcm$}%
$\zfc$ minus the Power Set axiom,  
%+ \lap{$\rV=\rL$},
and with the Collection rather than Replacement 
scheme,
and with the Wellorderability principle instead of 
the usual Axiom of Choice.
See \cite{gitPWS} about this. 

The theorem, sans the last 
claim (*), was established by 
Jensen~\cite{jenmin} by means of a forcing notion
$P\in\rL\yi P\sq\Lomi$,  
(called \rit{the Jensen minimal-$\id13$-forcing}, 
or simply \rit{Jensen's forcing}), 
and, without the minimality claim \ref{gt5},
also by Jensen and Solovay \cite{jsad}, by means 
of almost-disjoint forcing. 
The forcing notion $P$ is defined in \cite{jenmin} 
in $\rL$ in the form $P=\bigcup_{\al<\omi}P_\al$, 
where each $P_\al$ is a countable collection of 
perfect trees that has a certain relation to 
the lower levels $P_\ba$, $\ba<\al$, and at each 
step the extension $P_\al$ is basically the 
G\"odel-least of all suitable extensions. 
(See also \cite[28.A]{jechmill} for another but close
definition of $P$.)

It is well-known that there is a substantial family of
rather elementary forcing notions such that their
construction and basic forcing properties can be
reformulated, formalized {\em mutatis mutandis},
and proved on the basis of
$\zfcm$ and/or $\pad$.
The consequence of this is that whatever is forced
about the reals by such a forcing notion, can be
equiconsistently adjoined to $\zfcm$ and/or $\pad$,
with no need to use theories stronger in the sense of
relative consistensy.
This family of forcing notions includes Cohen, random,
Sacks, Silver and some other forcing notions.
Does it include the Jensen forcing $P$ defined in
\cite{jenmin}?

The original construction of $P$
(either by Jensen~\cite{jenmin} or via $\Diamond_{\omi}$
as in \cite[28.A]{jechmill}, does not substantiate the 
positive answer because
the proof of the key CCC property by either method
essentially depends on  
countable elementary submodels of $\rL_{\om_2}$,
which is way beyond $\zfcm$.
In this paper, we show how to circumwent this difficulty
by means of the method of killing only definable
antichains in the course of a Jensen-like transfinite
construction of the forcing notion.
Our changes are concentrated in
Definition~\ref{zfcpo} and Condition \ref{j3} in
Section~\ref{jf1+1}.

Then we define a
model with a minimal $\varPi^1_2$ singleton as a
class-forcing extension of a model of $\zfcm$ plus
$\rV=\rL$ plus ``all sets are countable''.
This will be the proof of Theorem~\ref{gt}.

\cha{Preliminaries} 
\las{2prel}

Let $\nse$ be the set of all strings (finite sequences) 
of natural numbers.
\index{strings!omegao@$\nse$}%   
\index{zzomega2omega@$\nse$}%   
Accordingly, $\bse\sneq\nse$ is the set of all dyadic strings  
(finite sequences of numbers $0,1$). 
\index{strings!omega2@$\bse$}%   
\index{zzomega2omega@$\bse$}%  
If $t\in\nse$ and $k<\om$ then 
$t\we k$ is the extension of $t$ by $k$ as the rightmost term. 
If $s,t\in\nse$ then $s\sq t$ means that $t$ extends $s$, while 
$s\su t$ means proper extension.%
\index{strings!extension, $s\su t$}%   

If $s\in\nse$ then $\lh s$ is the length of $s$,  
\index{zzlh@$\lh$}%
and $\om^n=\ens{s\in\nse}{\lh s=n}$ (strings of length $n$), 
and accordingly $2^n=\om^n\cap\bse=\ens{s\in\bse}{\lh s=n}$.
\index{strings!nomega@$\om^n$}%   
\index{zznomega@$\om^n$}%   
\index{strings!n2@$2^n$}%   
\index{zzn2@$2^n$}%   
 
A set $T\sq\nse$ is a \rit{tree} iff 
\index{tree}%
%it is an initial segment, that is, 
for any strings $s\su t$ in $\nse$, if $t\in T$ then $s\in T$.
Thus every non-empty tree $T\sq\nse$ contains the empty 
string $\La$. 
If $T\sq\nse$ is a tree and $s\in T$ then put 
\index{zzTrezs@$T\ret s$}%
\index{restrictionTress@restriction $T\ret s$}
$T\ret s=\ens{t\in T}{s\sq t\lor t\sq s}$; 
this is a tree as well.

Let $\pet$ be the set of all \rit{perfect} trees 
$\pu\ne T\sq \bse$. 
\index{zzpt@$\pet$}%
\imar{pet}%
Thus a non-empty tree $T\sq\bse$ belongs to $\pet$ iff 
it has no endpoints and no isolated branches. 
In this case, there is a largest string $s\in T$ such that 
$T=T\ret s$; it is denoted by $s=\roo T$   
\index{root of a tree}%
\index{zzrooT@$\roo T$}%
(the {\it root\/} of a perfect tree $T$). 
If $s=\roo T$ then $s$ is a \rit{branching node} of $T$, 
that is, 
\index{branching node}%
 $s\we 1\in T$ and $s\we 0\in T$.

Each perfect tree $T\in\pet$ defines a perfect set 
$$
[T]=\ens{a\in\dD}{\kaz n\,(a\res n\in T)}\sq\dD 
\index{zz[T]@$[T]$}
$$ 
of all paths through $T$; then accordingly 
$T=\ctr{[T]}$, where 
$$
\ctr X=\ens{a\res n}{a\in X\land n\in\om}\sq\bse 
\index{zztreeX@$\ctr X$}%
\quad\text{for any set}\quad
X\sq\omd.
$$
If $S\sq T$ are trees in $\pet$ and there is a finite set $A\sq T$ 
\index{tree!clopen}
such that $S=\bigcup_{s\in A}T\ret s$ then we say that $S$ is 
\rit{clopen} in $T$; 
then $[S]$ is a relatively clopen subset of $[T]$.
Trees clopen in $\bse$ itself will be called simply \rit{clopen}; 
thus clopen trees are those of the form 
$S= \bigcup_{s\in A}[s]$, where $A\sq\bse$ is a finite set and 
$[s]=\ens{t\in\bse}{s\sq t\lor t\sq s}$ for each $s\in\bse$.
\index{zz[s]@$[s]$}

A set $A\sq\pet$ is an \rit{antichain} iff $[T]\cap[S]=\pu$ 
\index{antichain}
(or equivalently, $S\cap T$ is finite) 
for all $S\ne T$ in $A$.
If $X\sq\pet$ then a set $D\sq X$ is: 
\bit
\item[$-$] \rit{dense in\/ $X$}, iff 
\index{set!dense}%
for every tree $T\in X$ there is a subtree $S\in D\yt S\sq T$; 

\item[$-$] \rit{open dense in\/ $X$}, iff 
\index{set!open dense}%
it is dense in $X$ and $T\in D$ holds
whenever $T\in X$, $S\in D$, and $T\sq S$;

\item[$-$] \rit{pre-dense in\/ $X$}, iff 
\index{set!pre-dense}%
the set $D^+=\ens{T\in X}{\sus S\in D\,(T\sq S)}$ is dense in $X$.
\eit
As usual, if $T\in\pet$, $D\sq\pet$, and there is a finite set  
$A\sq D$ such that $T\sq\bigcup A$
(or, equivalently, $[T]\sq\bigcup_{S\in A}[S]$)
then we write $T\sqf\bigcup D$,  
\index{zzsqf@$\sqf\bigcup$}%
\index{zzsqfa@$\sqfa\bigcup$}%
and if in addition $A$ is an antichain 
then we write $T\sqfa\bigcup D$.  

Thus perfect sets in the Cantor space $\dC=\omd$ 
\index{zzdd@$\dC$}%
are straightforwardly coded by perfect trees in $\pet$. 
It takes more effort to introduce a reasonable coding system 
for continuous functions $F:\dD\to\omm$. 
Let $\fpt$ (functional perfect trees) be the set of all sets
$c\sq \bse\ti\nse$ such that
\ben
\aenu
\itlb{fpt1} 
if $\ang{s,u}\in c$ then $\lh s=\lh u$; 

\itlb{fpt2}\msur
${c}$ is a tree, that is, if $\ang{s,u}\in{c}$ and 
$n<\lh s=\lh u$ then $\ang{s\res n,u\res n}\in{c}$; 
                                                            
\itlb{fpt3}\msur
$\dom{c}=\bse$, that is,
$\kaz s\in\bse\:\sus u\in\nse\:(\ang{s,u}\in{c})$; 

\itlb{fpt4}\msur 
${c}$ has no endpoints, that is, 
if $\ang{s,u}\in{c}$ and $\ell\in\ans{0,1}$ then there is 
$k<\om$ such that $\ang{s\we\ell,u\we k}\in{c}$; 

\itlb{fpt5} 
for every $m$ there exists $k\ge m$ such that  
if $s\in 2^k$ then there is a string $u_s\in\om^m$ 
satisfying  
$\kaz u\in\nse\,(\ang{s,u}\in{c}\imp u_s\sq u)$.
\een
If $F:\dD\to\omm$ is continuous then the set 
${c}=\cod F$, where 
$$
\cod F=\ens{\ang{a\res n,F(a)\res n}}{a\in\dD\land n<\om}
\quad\text{for any map}\quad
F:\dD\to\omm,
$$
belongs to $\fpt$ 
(condition \ref{fpt5} represents the uniform continuity of 
$F$ defined on a compact space), and $\fk{\cod{F}}=F$, where 
\imar{fk c}%
$$
{\fk{c}}=\ens{\ang{a,b}\in \dD\ti\omm}
{\kaz m\:\ang{a\res m,b\res m}\in{c})}
\quad\text{for every}\quad
{c}\in\fpt\,.
\index{zzfc@$\fk c$}%
$$
(a function coded by ${c}$). 
Conversely if ${c}\in\fpt$ then $\cod{\fk{c}}={c}$. 
 
\vyk{
Let $\fc c$ be the continuous function coded by $c\in \omm$,  
\index{zzfc@$\fc c$}
so that 1) if $c\in\omm$ then $\fc c:\dD\to\omm$ is continuous, and 
2) for any continuous $F:\dD\to\omm$ there is (not necessarily unique) 
$c\in\omm$ with $F=\fc c$.
}

The following is a well-known fact:

\ble
\lam{cORb}
If\/ $T\in\pet$ and\/ ${c}\in\fpt$ then 
either there is a string\/ 
$s\in T$ such that the restriction\/ 
${\fk{c}}\res [T\ret s]$ is 
a constant, or there is a subtree\/ 
$S\in\pet\yt S\sq T$, 
such that the restriction\/ 
${\fk{c}}\res [S]$ is a injection.\qed
\ele

\cha{Splitting systems of trees} 
\las{sst}

If $T\in\pet$ and $i=0,1$ then let 
$\spl Ti=T\ret{r\we i}$, 
\index{splitting!Ttoi@$\spl Ti$}%
\index{zzTtoi@$\spl Ti$}%
where $r=\roo{T}$; 
obviously $\spl Ti$ are trees in $\pet$ as well.
Define $\spl Ts$ for $s\in\bse$ by induction on $\lh s$ so that 
$\spl T\La=T$ and $\spl T{s\we i}=\spl{\spl Ts}i$.
\index{splitting!Ttos@$\spl Ts$}%
\index{zzTtos@$\spl Ts$}%

A \rit{splitting system}    
\index{splitting!system of trees}% 
is any indexed set $\sis{T_s}{s\in\bse}$ of trees $T_s\in\pet$
satisfying
\ben
\Aenu 
\itla{1spst}
if $s\in\bse$ and $i=0,1$,  then 
\imar{1spst}
$T_{s\we i}\sq \spl{T_s}i$.
\een
It easily follows from \ref{1spst} that 
\ben
\Aenu 
\atc
\itla{2spst}\msur
$s\sq s'\imp T_{s'}\sq T_s$, \ and 

\itla{3spst} 
if $n<\om$ and strings $s\ne t$ 
belong to $2^n$ then $[T_s]\cap[T_t]=\pu$.
\een

The next lemma is one of most known applications of the splitting.

\ble
\lam{simples}
If\/ $\sis{T_s}{s\in\bse}$ is a splitting system 
then\/ $T=\bigcap_n\bigcup_{s\in 2^n}T_s$ is a perfect 
subtree of\/ $T_\La$, 
and\/ $[T]= \bigcap_n\bigcup_{s\in 2^n}[T_s]$.
In addition, we have\/ 
$[T\ret s]=[T]\cap[T_s]$ for all\/ $s$.
\qed
\ele

We proceed to several slightly more complicated applications. 

\ble
\lam{sstL1}
If\/ $\ens{T^n}{n<\om}\sq\pet$ then there exists a sequence of 
trees\/ $S^n\in\pet$ such that\/ $S^n\sq T^n$ for all\/ $n$ and\/ 
$[S^k]\cap[S^n]=\pu$ whenever $k\ne n$. 
\ele
\bpf
Obviously if $T,T'\in\pet$ then there are perfect trees $S\sq T$
and $S'\sq T'$ such that $[S]\cap[S']=\pu$. 
This allows us to easily define a system 
$\sis{T_s(k)}{s\in\bse,\:k<\om}$ of trees $T_s(k)\in \pet$
such that
\ben
\nenu
\vyk{\itla{stL0}
if $k<\om$ then $T_\La(k)=T^k$;
\imar{stL0}
}
\itla{stL1}
if $k<\om$ then $\sis{T_s(k)}{s\in\bse}$ is 
a splitting system consisting of subtrees of $T^k$;

\itla{stL2}
if $k<\ell<n<\om$ and $s,t\in 2^n$ then 
$[T_s(k)]\cap[T_t(\ell)]=\pu$.
\een
(The inductive construction is arranged so that, at each step $n$, 
we define all trees $T_s(k)$ with $k<n$ and $s\in2^n$ and also 
all trees $T_s(n)$ with $\lh s\le n$.)
Now we simply put $S^k=\bigcap_n\bigcup_{s\in2^n}T_s(k)$ for all $k$.
\epf

\ble
\lam{sstC}
If\/ $\ens{T^n}{n<\om}\sq\pet$  
and\/ $F:\dD\to\omm$ is continuous, then there exist perfect 
trees\/ $S^n\sq T^n$ such that either\/ $F(a)\nin\bigcup_n[S^n]$ 
for all\/ $a\in [S^0]$, or\/ $F(a)=a$ for all\/ $a\in [S^0]$. 
\ele
\vyk{
\bpf
We assume that the \lap{or} case fails, that is, if  
$S^0\sq T^0$ is a perfect tree then $F(a)\ne a$ for at least one 
$a\in [S^0]$.
Then, moreover, 
\ben
\fenu
\itla{stl*}
if $k<\om$ and $S^0\yi S$ are perfect trees, $S^0\sq T^0$, 
\imar{utochnit detali}
then there exist perfect trees $U\sq S$ and $U^0\sq S^0$   
such that 
$F(a)\nin[U^0]\cup[U]$ for all $a\in [U^0]$.
\een
This allows us to define a system of perfect trees
$\sis{U_t}{t\in\bse}$, splitting in the above sense, 
and a sequence of perfect trees $S^k\;(k\ge1)$,
such that
\ben
\nenu
\itla{stL0'}\msur
$U_\La\sq T^0$ and $S^k\sq T^k$ for all $k\ge1$;
\imar{stL0'}

\itla{stL1'}
if $a\in[U_\La]$ then $F(a)\nin[U_\La]$;

\itla{stL2'}
if $k\ge1$ and $a\in\bigcup_{t\in2^k}[U_t]$ then $F(a)\nin[S^k]$.
\een
If this is done then it remains to define 
$S^0=\bigcap_k\bigcup_{t\in2^k}U_t$.

Now let us explain the inductive construction 
of the trees $T_s(k)$.\vom

\rit{Step 0}. 
Using \ref{stl*}, we pick a tree $U_\La\sq T^0$ such that 
$F(a)\nin [U_\La]$ for all $x\in[U_\La]$, so that \ref{stL1'} 
holds.\vom

\rit{Step 1}. 
Note that $2^1$ consists of two strings $\ang0$ and $\ang1$ of 
length 1. 
Put $S_{\ang\ell}=\spl{U_\La}\ell$ for $\ell=0,1$. 
Using \ref{stl*} 2 times, 
we obtain perfect trees $U_{\ang\ell}\sq S_{\ang\ell}$, $\ell=0,1$, 
and $S^1\sq T^1$
%and $T'_\La(n)\sq T_\La(n)\yt n\ge2$, 
such that $F(a)\nin [S^1]$ whenever 
$a\in[U_{\ang0}]\cup[U_{\ang0}]$.\vom

\rit{Step 2}. 
Note that $2^2$ consists of four dyadic strings of 
length 2, and each string in $2^2$ has the form $\tau\we \ell$, 
where $\tau\in2^1$ and $\ell=0,1$. 
If $\tau\in2^1$ and $\ell=0,1$ then put 
$S_{\tau\we\ell}=\spl{U_t}\ell$. 
Using \ref{stl*} 4 times consecutively, 
%(the number of all triples $\ang{s,t,u}$ of strings $s,t,u\in2^2$), 
we obtain perfect trees $U_t\sq S_t$ for $t\in2^2$, and $S^2\sq T^2$, 
such that $F(a)\nin [S^2]$ whenever 
$a\in \bigcup_{t\in2^2}[U_t]$.\vom

And so on. 
This construction obviously yields a system of trees as required.
\epf
}

\bpf
\snos
{This proof, much shorter and more transparent than 
our original proof, was suggested by the anonymous 
referee, and we thankfully follow their advice.}
Assume that $F(a_0)\ne a_0$ for some $a_0\in [T^0]$. 
By continuity
of $F$ there are a clopen subtree $S \sq T^0$ 
and a clopen neighbourhood $A$ of $F(a_0)$ 
such that $F([S])\sq A$ and $[S] \cap A =\pu$. 
Hence, $F(a) \nin [S]$ for all $a \in [S]$. 
The compact set $X = F([S])$ is either countable 
or has a perfect subset. 

If $X$ is countable then let $S^0 = S$ and 
for every $n\ge1$ let $S^{n}\sq T^{n}$ be an
arbitrary perfect tree such that 
$[S^{n}]\sq [T^{n}]\bez X$.

Assume that there is a perfect tree $T$ 
such that $[T]\sq X$. 
By Lemma \ref{sstL1}, there are trees
$U^n \in \pet$ such that $U^0 \sq T$, 
$U^{n+1}\sq T^{n+1}$, and [$U^k] \cap [U^n] = \pu$ 
whenever $k \ne n$.
Choose $S^0\in \pet$ such that 
$[S^0] \sq [S] \cap F^{-1}([U^0])$ 
and let $S^{n+1} = U^{n+1}$.
\epf

\vyk{

\cha{Jensen's construction: overview} 
\las{jco}

Beginning the proof of Theorem~\ref{gt},
we list essential properties 
of Jensen's forcing $\dJ\in\rL$:\vom

1) ${\dJ}$ consists of perfect trees $T\sq\bse$ 
(a subset of the Sackc forcing);\vom

2) ${\dJ}$ forces that there is a unique \dd {\dJ}generic real;\vom

3) \lap{being a \dd {\dJ}generic real} is a $\ip12$ property;\vom

4) ${\dJ}$ forces that the generic real is (nonconstructible and) 
minimal.\vom

\noi
Thus ${\dJ}$ forces a nonconstructible $\ip12$ real 
singleton $\ans a$  over 
$\rL$, whose only element is, 
therefore, a $\id13$ real in $\rL[a]$.

Jensen~\cite{jenmin} defined a forcing ${\dJ}$ in $\rL$ 
in the form 
${\dJ}=\bigcup_{\al<\omi}{\dJ}_\al$, 
where each ${\dJ}_\al$ is a countable 
collection of perfect trees $T\sq\bse$. 
The construction of the \dd\omi sequence of sets 
${\dJ}_\al$ is arranged 
so that each ${\dJ}_\al$ is generic, in a certain sense, 
over the least 
transitive model of a suitable fragment of $\ZFC$,  
containing the subsequence $\sis{{\dJ}_\ga}{\ga<\al}$. 
A striking corollary of such a genericity is that ${\dJ}$ 
forces that there 
is only one \dd {\dJ}generic real. 
Another corollary consists in the fact that for a real 
$x\in\dD$ being \dd {\dJ}generic is equivalent to 
$x\in\bigcap_{\al<\omi}\bigcup_{T\in {\dJ}_\al}[T]$. 
The construction can be managed so that the whole 
sequence $\sis{{\dJ}_\al}{\al<\omi}$ is $\id12$, or, 
more exactly, 
$\id\hc1$ in $\rL$. 
Altogether, it follows that if $a\in\dD$ is a 
\dd {\dJ}generic real then 
$\ans a\in\ip12$ in $\rL[a]$, that is, $a\in\id13$ in 
$\rL[a]$, which 
is obviously the lowest possible level for a 
nonconstructible real. 
The minimality of \dd {\dJ}generic reals follows from another 
property of ${\dJ}$: given a tree $S\in {\dJ}$ and a 
continuous $F:\dD\to\omm$, 
there is a tree $T\in {\dJ}\yd T\sq S$ (a stronger condition) 
such that 
$F\res[T]$ is either a bijection or a constant. 
This is a brief description of Jensen's proof of 
Theorem~\ref{gt} sans \ref{gt6}.
}

\cha{Jensen's sequences} 
\las{jf1+1}

In this section, \rit{we argue in\/ $\rL$.}

\vyk{
Recall that $\zfcm$ is the theory
\index{zzZFCm@$\zfcm$}%
\index{theory!ZFCm@$\zfcm$}%
$\zfc$ minus the Power Set axiom,  
%+ \lap{$\rV=\rL$},
and with the Collection rather than Replacement 
scheme,
and with the Wellorderability principle instead of 
the usual Axiom of choice.
(See \cite{gitPWS}.)
}

\bdf
[in $\rL$]
\lam{zfcpo}
Suppose that $\al<\omi$ and $\sis{X_\ba}{\ba<\al}$ 
is a sequence of sets in $\Lomi$. 
We let 
\imar{muu{}}
$\muu{\sis{X_\ba}{\ba<\al}}$
\index{zzmuJal@$\muu{\sis{X_\ba}{\ba<\al}}$}%
be the least ordinal $\mu\yi\al<\mu<\omi$,
such that:
\ben
\nenu
\itlb{mu1}
$\rL_\mu$ contains $\sis{X_\ba}{\ba<\al}$, 

\itlb{mu2}
$\rL_\mu$ 
contains the truth set $\ist{\rL_\al}$ which consists
of all closed \dd\in formulas with sets in $\rL_\al$
as parameters, true in $\rL_\al$,

\itlb{mu3}
$\rL_\mu$ models $\zfcm$ for bounded formulas,
plus ``all sets are countable''.\qed
\een
\eDf

\bdf
[in $\rL$]
\lam{j45}
If $\al<\omi$ then let $\ang{\nut \al,\nuf\al,\nu_\al}$ 
\index{zzTal@$\nut\al$}%
\index{zzcal@$\nuf\al$}%
be the \dd\al th element of the set $\pet\ti\fpt\ti\omi$ 
in the 
sense of the G\"odel canonical wellordering of $\rL$.
\edf

Thus for any $T\in\pet$ and $c\in\fpt$ there exist 
uncountably many 
indices $\al<\omi$ such that $T=T_\al$ and $c=c_\al$. 

For any ordinal $\la\le\omi$, we let $\js\la$ 
(Jensen's sequences of length $\la$) 
be the set of all sequences 
\index{zzJSla@$\js\la$}
$\sis{J_\al}{\al<\la}$ of length $\la$, of 
{\ubf countable}
sets $J_\al\sq\pet$, satisfying the following 
conditions \ref{j0} --- \ref{j5}. 
\ben
\cenu
%\def\theenumi{(J\arabic{enumi})}
%\def\labelenumi{\theenumi}
%\addtocounter{enumi}{-1}
\itlb{j0}\msur 
$J_0$ consists of all clopen trees $\pu\ne S\sq\bse$, 
including the full tree $\bse$ itself. 
%, and their finite unions (= all non-empty clopen sets in $\dD$).

\itlb{j1}
If $\al<\la$, $T\in J_\al$, 
and $S\sq T$ is a perfect tree clopen in $T$, 
then $S\in J_\al$. 

\itlb{j2} 
If $\al<\la$ and $S\in J_{<\al} =\bigcup_{\ba<\al} J_\ba$ 
then there is a tree $T\in J_\al\yt T\sq S$.

\itlb{j3}
%Suppose that $\al<\la$. 
%Let $\Mp\al=\minp{\sis{J_\ba}{\ba<\al}}$. 
%\index{zzM+al@$\Mp\al$}%
%\imar{Mp al}%
%We require that 
%
If $\al<\la$, 
$T\in J_\al$, $D\in \MM{\sis{J_\ba}{\ba<\al}}$, 
$D\sq J_{<\al}$, $D$ is open dense in  $J_{<\al}$,  
%or $D=J_\ba$ for some $\ba<\al$, 
then $T\sqfa\bigcup D.$ 

\itlb{j4}
%If $\al<\omi$, ${c}\in\fpt\cap \Mp\al$, and $S\in J_{<\al}$,
%then there is $T\in J_\al$ such that $T\sq S$ and
If $\al<\omi$, $c=\nuf\al$, and $S=\nut\al\in J_{<\al}$,
then there is $T\in J_\al$ \sth\ $T\sq S$ and:
\bde
\item[{\ubf either}]  
we have ${\fk{c}}(a)\nin\bigcup_{T'\in J_\al}[T']$ 
 for all $a\in [T]$,%\\[0.5ex] 
\item[{\ubf or}] 
we have ${\fk{c}}(a) = a$ for all $a\in [T]$.
\ede
 
\itla{j5}
If $\al<\omi$, $c=\nuf\al$, and $S=\nut\al\in J_{<\al}$,
then there exists $T\in J_\al$ such that $T\sq S$ and
the restricted function ${\fk{c}}\res [T]$ is either 
a bijection or a constant. 
\een

\noi
Let $\js{<\la}=\bigcup_{\al<\la} \js\al$. 
\index{zzJSlessl@$\js{<\la}$}

\ble
[in $\rL$]
\lam{jden}
Suppose that\/ $\ba<\la\le\omi$ and\/ 
$\sis{J_\al}{\al<\la}\in\js\la$. 
Then\/ $J_\ba$ is pre-dense in the set\/ 
$J_{<\la}=\bigcup_{\al<\la}J_\al$.
\ele
\bpf
First, $J_\ba$ is dense in $J_{<\ba+1}$ by \ref{j2}.
Now, by induction on $\la$, 
suppose that $J_\ba$ is pre-dense in $J_{<\la}$. 
To check that $J_\ba$ remains pre-dense in 
$J_{<\la+1}=J_{<\la}\cup J_\la$, 
consider any tree $T\in J_{\la}$. 
By definition $J_\ba\in \MM{\sis{J_\ba}{\ba<\la}}$, 
and hence we have $T\sqfa\bigcup J_\ba$ by \ref{j3}. 
(Note that the set 
$J_\ba^+=\ens{S\in J_{<\la}}{\sus S'\in J_\ba\,(S\sq S')}$ 
belongs to $\MM{\sis{J_\ba}{\ba<\la}}$ and is open dense.)
It follows that there exist a tree $S\in J_\ba$ and a string 
$s\in T$ such that $T\ret s\sq S$. 
Finally, $T'=T\ret s\in J_\la$ by \ref{j1},
so $T$ is compatible 
with $S\in J_\ba$, as required.
\epf

\ble
[in $\rL$]
\lam{ccc}
Assume that\/  
$\sis{J_\al}{\al <\omi}\in\js\omi$. 
Then the corresponding forcing notion\/ 
${J} =\bigcup_{\al <\omi} {J}_\al\in\rL$
satisfies CCC in\/ $\rL$ \poo\ all
antichains\/ $A\sq J$ definable in\/
$\Lomi$ with parameters. 
%Therefore the cardinals are preserved in\/ 
%\dd{J}generic extensions  of\/ $\rL$.
\ele
\bpf
Suppose that $A\sq {J}$ is a  
maximal \dd{J}antichain, that is, a pre-dense set and 
if $S\ne S'$ belong to $A$ then 
there is no tree $T\in {J}$, $T\sq S\cap S'$.
As $A$ is definable,
assume that $A=\ens{T\in\Lomi}{\Lomi\mo\vpi(p,T)}$,
where $p\in\Lomi$ is a parameter and $\vpi$ any
\dd\in formula.

There exists a limit ordinal $\al$ such that
$p\in\rL_\al$, the set 
$J_{<\al}=\bigcap_{\ga<\al}J_\ga$ satisfies 
$J_{<\al}=J\cap\rL_\al$,
%$A\cap J_{<\al}=A\cap\rL_\al$,  
the set
$A_{<\al}=A\cap J_{<\al}$ is a maximal antichain,
therefore pre-dense in $J_{<\al}$, and finally
$\rL_\al$ is elementarily equivalent to $\Lomi$
\poo\ $\vpi$, so that overall we have:
$A_{<\al}=\ens{T\in\rL_\al}{\rL_\al\mo\vpi(p,T)}$.

Let $\mu=\muu{\sis{J_\ga}{\ga <\al}}$.
We assert that $A_{<\al}\in \rL_\mu$.
Indeed, by definition the truth set
$\btau=\ist{\rL_\al}$ belongs to $\rL_\mu$.
On the other hand,
$A_{<\al}=\ens{T}{\vpi(p,T)\in\btau}$
by the above.
It follows that $A_{<\al}\in \rL_\mu$ since $\rL_\mu$
models $\zfcm$ for bounded formulas.

Now it suffices to prove that $A=A_{<\al}$. 
Suppose towards the contrary that 
$T\in A\bez A_{<\al}=A\bez {J}_{<\al}$. 
Then $T$ is compatible with some $T'\in {J}_\al$ 
by Lemma~\ref{jden}, that is, there is 
a tree $T''\in {J}$, $T''\sq T'\cup T$.

On the other hand, it follows from \ref{j3}
that $T'\sqfa A_{<\al}$. 
Then $T''\sqfa A_{<\al}$ as well, and hence there exist 
$s\in T''$ and 
$S\in A_{<\al}$ such that the tree $U=T''\ret s$ satisfies 
$U\sq T''\cap S$, therefore $U\sq T\cap S$. 
However $U\in {J}$ by \ref{j1}, and $S\in A_{<\al}$ 
but $T\in A\bez A_{<\al}$, 
contrary to the assumption that $A$ is a 
\dd {J}antichain.
\epf 

The following 
rather obvious lemma demonstrates that the top level of a 
Jensen sequence of successor length can be freely enlarged 
by adding smaller trees, with only care of the 
property \ref{j1}.

\ble
[in $\rL$]
\lam{l4-}
Suppose that\/ $\la=\xi+1<\omi$ and\/ 
$\fj=\sis{J_\al}{\al<\la}\in\js\la$, so that\/ $J_\xi$ is 
the last set in\/ $\fj$. 
Assume that\/ $S\sq T$ are trees in\/ $\pet$ and\/ $T\in J_\xi$. 
Let\/ $J'_\xi$ consist of all trees in\/ $J_\xi$ and all trees\/ 
$S'\in\pet\yt S'\sq S$, clopen in\/ $S$. 
Then the sequence\/ 
$\fja=\sis{J_\al}{\al<\xi}\cup\ans{\ang{\xi,J'_\xi}}$ 
still belongs to\/ $\js{\la}$.\qed
\ele

\cha{Extension of Jensen's sequences}
\las{jext}

Now we prove a theorem 
which shows that Jensen's sequences of any countable 
length are extendable to longer sequences in $\rL$. 

\bte
[in $\rL$]
\lam{l4}
Suppose that\/ $\la<\omi$. 
Then any sequence\/ 
$\fj=\sis{J_\al}{\al<\la}\in\js\la$ has 
an extension\/ 
$\fja=\sis{J_\al}{\al\le\la}\in\js{\la+1}$.
\ete
\bpf
\rit{We argue in\/ $\rL$.}
Basically, we have to appropriately define the 
\rit{top level} $J_\la$ ($\la>0$) of the extended sequence. 
The definition goes on in four steps. 

{\ubf Step 1}: 
we define a provisional set $J_\la$ satisfying only 
requirements \ref{j2}, \ref{j3}. 
Put $\Mp\la=\MM{\sis{J_\ba}{\ba<\la}}$.
Fix an arbitrary enumeration $\ens{D_n}{n<\om}$ of 
all sets $D\in \Mp\la$, $D\sq J_{<\la}$, open  
dense in $J_{<\la}$, and an arbitrary enumeration
$J_{<\la}=\ens{S^k}{k<\om}$. 
For any $k$, there is a system 
$\sis{T_s(k)}{s\in\bse}$ of trees $T_s(k)\in J_{<\la}$
satisfying the following conditions 
\ref{spis1} -- \ref{spisL}:
\ben
\renu
\itla{spis1}
if $S=S^k\in J_{<\la}$ then $T_\La(k)\sq S$; 

\itla{spis2}\msur
$\sis{T_s(k)}{s\in\bse}$ 
is a splitting system in the sense of \nrf{sst}; 

\itla{spis3}
\label{spisL}
if $n\ge1$ and  $s\in 2^n$ then $T_s(k)\in D_{n}$.
\een
Indeed if some $T_s(k)\in J_{<\la}$ 
is already defined and $n=\lh s$, 
then the trees $U_0=T_s(k)\ret0$ 
and $U_1=T_s(k)\ret1$ belong to 
$J_{<\la}$ either, and hence there are trees 
$T_{s\we0}\sq U_0$ and 
$T_{s\we1}\sq U_1$ in $J_{<\la}$, 
which belong to $D_{n+1}$. 

It remains to define $J_\la=\ens{T^k}{k<\om}$, where 
$T^k=\bigcap_n\bigcup_{s\in2^n}{T_s(k)}$. 

{\ubf Step 2}. 
We shrink the trees $T^k$ obtained at Step 1 in order to 
satisfy requirement \ref{j4}. 
Note that if $J_\la=\ens{T^k}{k<\om}$ 
satisfies \ref{j2} and \ref{j3} 
and $U^k\sq T^k$ is a perfect tree 
for each $k$ then the new set 
$J_\la=\ens{U^k}{k<\om}$ still 
satisfies \ref{j2} and \ref{j3}. 

Now suppose that $c=\nuf\la$ and 
$S=\nut\la\in J_{<\la}$, as in 
\ref{j4}. 
(If $S\nin J_{<\al}$ then we skip this step.) 
We may assume that the enumeration 
$\sis{T^k}{k<\om}$ is chosen 
so that $T^0\sq S$. 
Let $G=\fk c$ (a continuous map $\dD\to\omm$).
By Corollary~\ref{sstC}, 
there exist perfect trees $U^n\sq T^n$ 
such that either $G(a)\nin\bigcup_n[U^n]$ 
for all\/ $a\in [U^0]$, 
or\/ $G(a)=a$ for all\/ $a\in [U^0]$. 
The new set $J_\la=\ens{U^k}{k<\om}$  still satisfies 
\ref{j2} and \ref{j3}.

{\ubf Step 3}.  
We shrink the trees $U^k\in\pet$ obtained at Step 2 
in order to 
satisfy \ref{j5}. 
This is similar to Step 2, with the only difference 
that we apply 
Lemma~\ref{cORb} instead of Corollary~\ref{sstC}. 

{\ubf Step 4}.  
If $V^k\in\pet$ is one of the trees in $J_\la$ obtained 
at Step 3 
then we adjoin all trees $\pu\ne S\sq V^k$ clopen in $V^k$, 
in order to 
satisfy \ref{j1}.
\epf

\cha{Definable Jensen's sequence}
\las{djs}

Each of the conditions \ref{j3}, \ref{j4}, \ref{j5} 
(Section~\ref{jf1+1})  
will have its own role. 
Namely, \ref{j3} implies ccc and continuous reading 
of names (Lemma~\ref{l2}), 
\ref{j4} is responsible for the generic uniqueness 
of $a_G$ as in Lemma~\ref{l3}, 
while \ref{j5} yields the minimality. 
However, to obtain the required type of definability of 
\dd{\dJ}generic reals in the extensions, we need to take care of   
appropriate definability of a Jensen's sequence 
%(see Definition~\ref{fixJ}) 
in $\rL$. 

\bdf
\lam{1hc}
$\hc$ is the collection of all
{\it hereditarily countable\/} sets. 
\index{zzHC@$\hc$}
Note that $\HC=\Lomi$ under $\rV=\rL$. 
\bde
\item[$\hcs n$] = all sets 
\index{zzSihcn@$\hcs n$}%
\index{zzPihcn@$\hcp n$}%
\index{zzDehcn@$\hcd n$}%
\index{definability class!Sihcn@$\hcs n$}%
\index{definability class!Pihcn@$\hcp n$}%
\index{definability class!Dehcn@$\hcd n$}%
$X\sq\hc$, definable in $\hc$ by a parameter-free $\Sg_n$ formula. 
\vyk{
\item[$\bs n$] = all sets $X\sq\hc$  
\index{zzSihcbn@$\bs n$}%
\index{zzPihcbn@$\bp n$}%
\index{zzDehcbn@$\bd n$}%
\index{definability class!Sihcbn@$\bs n$}%
\index{definability class!Pihcbn@$\bp n$}%
\index{definability class!Dehcbn@$\bd n$}%
definable in $\hc$ by a $\Sg_n$ formula 
with sets in $\hc$ as 
parameters.
}%
\ede 
Collections $\hcp n\yt \hcd n$
%, {\it etc.\/}  
\index{zzSihcnx@$\hcs n(x)$}%
\index{zzPihcnx@$\hcp n(x)$}%
\index{zzDehcnx@$\hcd n(x)$}%
\index{definability class!Sihcnx@$\hcs n(x)$}%
\index{definability class!Pihcnx@$\hcp n(x)$}%
\index{definability class!Dehcnx@$\hcd n(x)$}%
are defined similarly. 
%
%Something like $\hcs n(x)\yt x\in\hc$, 
%means that only $x$ is admitted as a parameter.%
Essentially $\hcs n\yt\hcp n\yt\hcd n$ is the same as 
$\is1{n+1}\yt\ip1{n+1}\yt\id1{n+1}$ for sets of reals, 
modulo any appropriate coding.
%, and the same for boldface classes.
%
\edf

\bpro
[in $\rL$]
\lam{cl3}
The set\/ 
$\ens{\ang{\al,\fj}}{\al<\omi\land \fj\in\js\al}$ 
is\/ $\id\Lomi1$.
\epro
\bpf
Straightforward analysis of the definitions in
Section~\ref{jf1+1}.
\epf

\vyk{
Suppose that $J$ is a sequence (of any kind) 
of length $\la<\omi$, 
and $\vt<\omi$ is such that the set $\rL_\vt$ contains $J$, 
is a model of 
$\zfcm$, and for every $\al<\la$ the model 
$\Mp\al=\minm{\sis{J_\ba}{\ba<\al}}$ 
%defined in \ref{j3} of \nrf{jf1,} 
also belongs to $\rL_\vt$. 
Then the property of $J$ being a Jensen sequence is   
absolute for $\rL_\vt$. 
This yields a $\is{}1$ definition for the statement
\lap{$J$ is a Jensen sequence}
in the form: there is such-and-such ordinal $\vt$ such that 
$J\in\rL_\vt$ and 
\lap{$J$ is a Jensen sequence} holds in $\rL_\vt$.
}

\bcor
[in $\rL$]
\lam{e2}
There exists a $\id\Lomi1$ sequence\/ 
$\fj=\sis{\dJ_\al}{\al <\omi}\in\js\omi$. 
\ecor
\bpf
For every $\al,$ we let $\dJ_\al$ to be the least set, 
in the sense of G\"odel's $\id{}1$ wellordering of $\Lomi$, 
such that $\sis{\dJ_{\ba}}{\ba\le\al}\in\js{\al+1}$. 
\epf

\cha{Adding one Jensen real} 
\las{jf2}

Here we prove Theorem~\ref{gt} without claim (*). 

\bdf
\lam{n=2J}
Fix a sequence $\xj=\sis{\dJ_\al}{\al<\omil}\in\rL$, 
such that it is true in $\rL$ that\vom 

1) 
$\sis{\dJ_\al}{\al<\omi}\in\js{\omi}$, and\vom 

2)
$\sis{\dJ_\al}{\al<\omi}$ is a $\id\Lomi1$ sequence.\vhm 

\noi
(We refer to Corollary~\ref{e2}.)
Put $\dJ =\bigcup_{\al<\omi} \dJ_\al$. 
\edf

Consider such a set $\dJ\sq\pet$ as a forcing notion over 
$\Lomi$. 
It is ordered so that 
$S\sq T$ means that $S$ is stronger as a forcing condition. 
Thus $\dJ$
%, \rit{the Jensen forcing} of 
%\cite{jenmin} (see also \cite[28.A]{jechmill}), 
consists of (some, not all) perfect trees. 
Forcing notions of this type are sometimes called 
\rit{arboreal}.
 
\ble
%[in $\rL$]
\lam{cl1+}
If\/ $G\sq\dJ$ is a\/ $\dJ$-generic set over\/ $\Lomi$,
then the 
intersection\/ $\bigcap_{T\in G}[T]$ is a singleton\/ 
$\ans{a_G}\yt a_G\in\dD$, and\/ $G =\ens{T\in \dJ}{a_G\in T}$, 
hence\/ $\Lomi[a_G] = \Lomi[G]$. 
\ele

\bpf
Make use of \ref{j1}. 
\epf 

Reals $a_G$, $G\sq {\dJ}$ being a \dd {\dJ}generic set over 
$\Lomi$, are called {\it\dd {\dJ}generic over\/} $\Lomi$.

The following lemma provides a useful method of 
representation for 
reals in \dd {\dJ}generic extensions. 

\ble
[continuous reading of names]
\lam{l2}
Suppose that\/ $G\sq {\dJ}$ is\/ \dd {\dJ}generic over\/
$\Lomi$.
Let\/ $x\in \Lomi[G]\cap\omm$. 
There exists\/ ${c}\in\Lomi\cap \fpt$ such 
that\/ $x = \fk{c}(a_G)$. 
\ele
\bpf
The proof is based mainly on \ref{j3}. 
Let, indeed, $\dox$ be a 
name for $x$ in the forcing language, so that
every $T\in {\dJ}$ forces $\dox\in \omm$, and
$$
x(n)=l\;\eqv\; \sus T\in G\;(T\wfo \dox(\don)=\dol). 
$$
Let $T_0\in {\dJ}.$ 
We define in $\Lomi$, the ground model, 
$$
D_{nl} =\ens{T\in {\dJ}}{T\wfo \dox(\don)=\dol},
\quad\text{and}\quad 
D_n ={\textstyle\bigcup_{l\in\om}} D_{nl}. 
$$
All sets $D_n$ are dence in ${\dJ}.$ 
Arguing as in the proof 
of Lemma~\ref{ccc}, we obtain an 
ordinal $\al<\omi$ such that
$T_0\in {\dJ}_{<\al}=\bigcup_{\ga<\al}{\dJ}_{\ga}$,  
and, for any $n,$ the set $D_n(\al)=D_n\cap {\dJ}_{<\al}$ 
belongs to $\rL_\mu$,
where $\mu=\muu{\sis{\dJ_\ga}{\ga <\al}}$, 
and is dense in ${\dJ}_{<\al}$. 
By \ref{j2}, there exists $T\in {\dJ}_\al\yd T\sq T_0$. 
By \ref{j3} we have $T\sqfd\bigcup D_n(\al)$ for every $n$, 
so that there are finite sets $D'_n\sq D_n(\al)$ such that
$T\sq\bigcup D'_n$ and if $S\ne S'$ belong to the same set 
$D'_n$ then $[S]\cap[S']=\pu$. 

We put $D'_{nl}=D'_n\cap D_{nl}$. 
Then for any $n$ there is a finite number of values of $l$ 
such that $D'_{nl}\ne\pu$. 
Thus one can define in $\Lomi$ a continuous 
function $F':[T]\to \omm$ as follows: 
$F'(x)(n)=l$ iff $x\in [T]$ for some $T\in D'_{nl}$. 
Let $F:\dD\to\omm$ be a continuous extension of 
$F'$; $F=\fk{c}$ for some ${c}\in \fpt\cap\Lomi$. 
Then $T$ forces 
$\dox=\fk {\doc}(\doa)$, where $\doa$ is the canonical  
name for $a_G$. 
\epf

%Assume that \ref{j4} holds (in $\rL$) --- in addition to 
%Definition~\ref{fixJ}.

\ble
%[uniqueness] 
\lam{l3}
If\/ $G\sq {\dJ}$ is a\/ ${\dJ}$-generic set over\/
$\Lomi$ then\/
$a = a_G$ is the only element of\/ 
$\bigcap_{\al<\omi}\bigcup_{T\in {\dJ}_\al}[T]$ in\/
$\Lomi[G]$. 
Moreover\/ $a_G$ is minimal over\/ $\Lomi$.
\ele
\bpf
If $\al<\omi$ then 
the real $a=a_G$ actually belongs to 
$\bigcup_{T\in {\dJ}_\al}[T]$ 
since all sets ${\dJ}_\al$ are pre-dense by Lemma~\ref{jden}. 
To prove the opposite direction, consider any $S\in {\dJ}$ 
and $b\in 2^\om\cap\Lomi[G]$.
By Lemma~\ref{l2}, there exists ${c}\in\Lomi\cap \fpt$ such 
that $b = \fk{c}(a_G)$. 
There is an ordinal $\al<\omi$ in $\Lomi$ such that $T=\nut\al$ 
and $c=\nuf\al$. 
Let $T\in {\dJ}_\al$ witness \ref{j4}. 
In the \lap{either} case of \ref{j4}, $T$ obviously forces that 
$\fk{c}(a_G)\not\in\bigcup_{T'\in {\dJ}_\al} [T']$, while in 
the \lap{or} case $T$ forces $\fk{c}(a_G) = a_G$.

To prove the minimality consider any real 
$b\in2^\om\cap\Lomi[a_G]$. 
By Lemma~\ref{l2} we have $b=\fk{c}(a_G)$, 
where ${c}\in\fpt\cap\Lomi$. 
It follows from \ref{j5} that there exists 
$T\in G$ such that 
$\fc c\res[T]$ is either a bijection or a constant. 
If $\fk{c}\res[T]$ is a bijection then 
$a_G\in\Lomi[x]$ by means of the inverse map. 
If $\fk{c}\res[T]$ is a constant $z$, 
say $\fk{c}(x)=z$ for all $x\in[T]$ in $\Lomi$, 
then obviously $b=\fk{c}(a_G)=z\in\Lomi$.
\epf

Now consider any set $G\sq {\dJ}$ \dd {\dJ}generic over
$\Lomi$. 
Lemma~\ref{l3} implies 
$\ans{a_G}\in\ip\hc1,$ hence $\in\ip12,$ in
$\Lomi[G]=\Lomi[a_G]$. 
Thus $a_G\in\id13$ in $\Lomi[a_G]$, as required. 

This completes the proof of Theorem~\ref{gt} sans
claim (*) of the theorem.

\cha{Down to $\zfcm$}
\las{zfcm}

Now let's argue in the theory
$$
\text{
$\zfcm$ plus $\rV=\rL$ plus
``all sets are countable'',
}
\eqno(\dag)%
$$%
Its universe can be identified with $\Lomi$.
The above construction is basically
relativized to $\Lomi$, so that it can be executed in
the universe of $(\dag)$, which we denote by
$\Lomid$ for the sake of convenience.

Then $\dJ=\bigcup_{\al}\dJ_\al$ is a definable class,
more exactly $\id{}1$, and a class-forcing notion, CCC
\poo\ all definable (with parameters) class-antichains.
It is known (see \eg\ \cite{AG}) that this suffices to
develop forcing engine to the extent of making valid
in this setting all suitable results valid in the
context of \dd\dJ forcing over $\Lomi$.

We conclude that \dd\dJ generic extensions of $\Lomid$
prove Theorem~\ref{gt}
(including claim (*) of the theorem).\vtm

\qeD{Theorem~\ref{gt}}

\vyk{

\cha{Warmup: minimal $\id1n$ reals}
\las{jf:n}

To solve the case of an arbitrary index $n\ge3$ in 
Theorem~\ref{mt2}, we employ one more idea. 
Jensen's \dd\omi sequence $\sis{{\dJ}_\al}{\al<\omi}$ 
as in \ref{n=2J} 
can be seen as an 
\dd\omi branch of class $\id\hc1$ in the set $\jti$ of all 
countable (transfinite) sequences satisfying conditions  
\ref{j1} -- \ref{j5} above. 
%, that can be explicitly 
%extracted from the arguments that work for $n=2$. 
Now, fix $n\ge3$ and consider another \dd\omi branch 
$\sis{{\dJ}_\al}{\al<\omi}$ 
in the same set $\jti$, which is \rit{complete}, or \rit{generic}, 
in the sense that it intersects any $\bs{n-2}$ set 
$D\sq\jti$ dense in $\jti$.\snos
{The actual requirement will be slightly different.} 
%Recall that $\bs{n-2}$, or sometimes $\is{}{n-2}(\hc)$, 
%means the definability in $\hc$ by a 
%$\is{}{n-2}$ formula with arbitrary parameters in $\hc$. 
%
Such a sequence cannot be $\id\hc1$ in $\rL$,
but it is easy to get it in $\id\hc{n-1}$. 
It follows that in this case any \dd {\dJ}generic real $a$ is $\id1{n+1}$ 
in $\rL[a]$. 
Moreover, it follows from the completeness that ${\dJ}$ is an 
\lap{elementary subforcing} of the Sacks forcing $\pet$ in the sense 
that both force the same $\is1n$ sentences.
The Sacks forcing, being homogeneous, decides any sentence with names 
for \lap{old} sets as parameters. 
So it follows that ${\dJ}$ decides any $\is1n$ sentence, as required.

Now let us present this construction in detail.

\vyk{ 

It happens that Jensen's construction of ${\dJ}$ 
can be carried oud so that, 
for a given $n\ge3$, ${\dJ}$ is elementarily equivalent to 
$\per$ \poo\ the 
forcing of $\is1n$ formulas, while the above 
properties 1), 2), 3), 4) 
are preserved, with only 3) 
being accordingly weakened to $\ip1n$.

We present here both Jensen's construction in Section~\ref{1jf} 
and the 
modification in Section~\ref{jf2}. 
}

\cha{Complete sequences and sets of trees}
\las{thi} 

Coming back to the general case of 
Theorem~\ref{mt2}, we 
introduce a couple of important definitions.

\bdf
\lam{solv}
Suppose that $P=\stk P\pce$ is a partially ordered set.
For any set $D\sq P$ let $\sol D{}=\sol DP$ be the set of all 
$p\in P$ such 
\index{solves a set}%
\index{zzDsolvP@$\sol DP$}%
that either $p\in D$ or there is no $q\in D\yt q\pce p$. 
\vyk{
\snos
{Using Definition~\ref{solv}, it is usually the case that 
$p\pce q$ means that either $p$ is a stronger forcing condition 
in $P$ or that $p,q$ 
are sequences of some kind and $p$ extends $q$.} 
}%
\edf

\bdf
[in $\rL$]
\lam{dfnc}
Suppose that $n\ge 3$.
A sequence $\sis{J_\al}{\al<\omi}\in\js\omi$ is 
\dd n{\it complete} if for any $\bs{n-2}$ set $D\sq \jti$ 
there is $\ga <\omi$ such that 
$\sis{J_\al}{\al<\ga}\in \sol D{}$ --- 
meaning that either $\sis{J_\al}{\al<\ga}\in D$ 
or there is no sequence in 
$D$ extending $\sis{J_\al}{\al <\ga}$. 

A set of perfect trees $J\sq\pet$ is \dd n{\it complete} if for 
%any $T\in J$ and 
any $\bs{n-2}$ set $W\sq \pet$, the set 
$\sol W{}\cap J=
\ens{S\in J}{S\in W\lor \neg\:\sus T\in W\,(T\sq S)}$ 
is dense in $J$. 
%
%there is a tree $S\in J\yt S\sq T$, such that either $S\in W$, 
%or there is no tree $S'\in W$ with $S'\sq S$.
\edf

Thus $n$-completeness is a property of ``generic" nature, where 
genericity is related to a family of sets distinguished by a 
definability property. 

\ble
[in $\rL$]
\lam{dd}
Suppose that $n\ge 3$.
If a sequence\/ $\sis{J_\al}{\al<\omi}\in\js\omi$ is\/ 
\dd ncomplete then the set\/ $J=\bigcup_{\al<\omi}J_\al$ is\/ 
\dd ncomplete.
\ele
\bpf\snos
{We are thankful to the anonymous referee for correcting 
this proof.}
Suppose that $W\sq\pet$ is a $\bs{n-2}$ set, and 
$S\in J$, that is, $S\in J_\vt$ for some $\vt<\omi$. 
We prove that there is $T \in \sol W{}\cap J$ 
such that $T \sq S$.
The set $D$ of all sequences 
$\sis{J'_\al}{\al<\la}\in\jti\yt\la<\omi$, 
such that there exists 
$T\in\bigcup_{\al<\la}J'_\al\cap W\yt T\sq S$, is $\bs{n-2}$.
It follows that $\sis{J_\al}{\al<\la}\in\sol D{}$ 
for some $\la<\omi$, 
\ie, either $\sis{J_\al}{\al<\la}\in D$, 
or there is no sequence 
in $D$ that extends $\sis{J_\al}{\al<\la}$. 

If $\sis{J_\al}{\al<\la}\in D$ then by definition there 
exists a tree $T\in J\cap W$ with $T\sq S$, and 
$T$ is as required.

Now suppose that $\sis{J_\al}{\al<\la}$ 
is not extendable to a sequence in $D$, and 
denote $\xi=\max\ans{\la,\vt+1}$. 
Then the extended sequence 
$\fj=\sis{J_\al}{\al\le\xi}$ 
is not extendable to a sequence in $D$ because 
$\sis{J_\al}{\al<\la}$ is not extendable. 
By \ref{j2} there is a tree $T\in J_\xi$, $T\sq S$. 
We claim that $T\in\sol W{}$. 

% that is, in this case, there is 
%no tree $T'\in W\yt T'\sq S$ --- which obviously implies 
%the lemma. 

Suppose towards the contrary that $T\nin W$ and 
there is $T'\in W$ such that $T'\sq T$. 
Then by Lemma~\ref{l4-} there is 
a set $J'_\xi\sq\pet$ containing $T'$ and such that 
$\fj'=\fj\cup\ans{\ang{\xi,J'_\xi}}$ is 
still a sequence in $\js{\xi+1}$ extending $\fj$, and 
$\fj'\in D$ by the choice of $T'$.
But this contradicts the 
non-extendability of $\fj$, and therefore $T\in\sol W{}$.
\epf

\ble
[in $\rL$]
\lam{l5}
If\/ $n\ge3$ then there exists an\/ \dd ncomplete\/ 
sequence\/ $\sis{\dJ_\al}{\al<\omi}\in\js\omi$ of class\/ 
$\id\hc{n-1}$. 
\ele
\bpf
Let $U\sq\hc\ti\hc$ be a universal $\is\hc{n-2}$ set. 
That is, $U$ itself is $\is\hc{n-2}$, and if 
$X\sq\hc$ is a (boldface) $\bs{n-2}$ set then there 
is a parameter $p\in\hc$ such that 
$X=U_p:=\ens{x\in\hc}{\ang{p,x}\in U}$. 
As we argue in $\rL$, for any $\al<\omi$ 
let $p_\al$ be the $\al$th element of $\hc=\rL_{\omi}$ 
in the sense of G\"odel's $\id{\hc}1$ wellordering 
of $\hc=\rL_{\omi}$. 
Then $\hc=\ens{p_\al}{\al<\omi}$ and the sequence 
$\sis{p_\al}{\al<\omi}$ is $\is\hc1$. 

To prove the lemma, we define a strictly 
\dd{\subset}increasing sequence 
$\sis{{\boldsymbol j}[\al]}{\al<\omi}$ 
of sequences ${\boldsymbol j}[\al]\in\jti$ as follows. 
Let ${\boldsymbol j}[0]$ be the empty sequence. 

Let 
${\boldsymbol j}[\la]=\bigcup_{\al<\la}{\boldsymbol j}[\al]$ 
whenever $\la<\omi$ is limit. 

For every $\al,$ if ${\boldsymbol j}[\al]\in\jti$ 
is defined, then let ${\boldsymbol j}[\al+1]$ to be  
the G\"odel-least sequence ${\boldsymbol j}\in\jti$, 
such that ${\boldsymbol j}[\al]\sq{\boldsymbol j}$ 
and ${\boldsymbol j}\in\sol{U_{p_\al}}{}$. 

The limit sequence 
$\sis{\dJ_\al}{\al<\omi}=
\bigcup_{\al<\omi}{\boldsymbol j}[\al]\in\js\omi$ 
is $n$-complete by construction, and an easy estimation, 
based on the assumption that $U$ is $\is\hc{n-2}$, 
shows that it belongs to $\id\hc{n-1}$. 
\epf

The next theorem is the main step in the proof of
Theorem~\ref{mt2}.
Its proof will be accomplished in \nrf{pt5}.
 
\bte
[in $\rL$]
\lam{mt2+}
Suppose that\/ $n\ge3$, $\sis{{\dJ}_\al}{\al<\omi}\in\js\omi$ is 
an\/ \dd ncomplete sequence of class\/ $\id\hc{n-1}$ 
(by Lemma~\ref{l5}), 
and\/ ${\dJ} =\bigcup_{\al <\omi } {\dJ}_\al$.
Then\/ \dd {\dJ}generic extensions of\/ $\rL$ provide  the proof of 
Theorem~\ref{mt2}.
\ete 

A few remarks before the proof starts.

It follows from Lemma~\ref{l3} that if a set $G\sq {\dJ}$ is 
\dd {\dJ}generic over $\rL$ then the corresponding real $a_G$ is 
minimal. 
It also follows from the same lemma and the fact that the 
sequence $\sis{{\dJ}_\al}{\al<\omi}\in\js\omi$ is $\hcd{n-1}$ in 
$\rL$ that the singleton $\ans{a_G}$ is $\ip1n$ and hence $a_G$ 
itself is $\id1{n+1}$ in $\rL[G]$.
It is a more difficult problem to 
prove the remaining claim, that is, that any 
$\is1n$ set $x\sq\om$ in $\rL[G]$ is constructible. 
We'll establish this fact in the remainder; the result 
will be based on the \dd ncompleteness property and on some 
intermediate claims.

\cha{Definability of the Sacks forcing}
\las{dsf} 

%Prior to the general case $n\ge 3$ of Theorem~\ref{mt2}, 
%we have to take a pause for the 
Our next goal is to estimate the definability of the 
Sacks forcing relation, restricted to formulas of a certain 
ramified version of the 2nd order Peano language.  

\bdf
\lam{sfor}
Let $\ya$ be the ordinary language of the 2nd order 
Peano arithmetic, 
\index{formula!L@$\ya$}%
with variables of type 1 for \rit{functions} in $\om^\om$.
Extend this language so that 
some type 1 variables 
can be substituted by symbols of the form $\ovc\yd c\in\fpt$, 
\index{name!c@$\ovc$}%
\index{zzchat@$\ovc$}
and each $\ovc$ is viewed as a name for $\fk{c}(a)$, where $a$ means 
a generic real of any kind. 
(Recall that $\fk{c}:\dD\to\omm$ is a continuous map coded by 
${c}\in\fpt$.)
Let $\sya$ be the extended language; the index s is from Sacks. 
\index{formula!Sind@$\sya$}%
Accordingly, $\cis1n$ and $\cip1n$ will denote the standard classes 
of formulas of the extended language $\sya$.
\index{formula!Si1ks@$\cis1k$}%
\index{formula!Pi1ks@$\cip1k$}%
\index{zzsSi1ks@$\cis1k$}%
\index{zzsPi1ks@$\cip1k$}%

If $a\in\dD$ and $\vpi$ is a formula of $\sya$  
then $\vpi[a]$ is the result of substitution of $\fk c(a)$ 
\index{zzfiG@$\vpi[a]$}%
for any name $\ovc$ in $\vpi$; 
$\vpi[a]$ is a formula of $\ya$ with real parameters. 
\edf

%The following definition introduces a somewhat modified 
%\dd\pet forcing notion $\bfo$ for formulas in $\sya$. 

\bdf
\lam{sforc}
Let $\for$ be the Sacks forcing relation 
(that is, $\pet$ is the forcing notion). 
% over $\rL$.   
\index{relation!for@$\for$}%
\index{zzfor@$\for$}%
Define an auxiliary relation of \lap{strong} forcing $\bfo$,    
\index{relation!forcS@$\bfo$}%
\index{zzforcS@$\bfo$}%
restricted to $\cis1k$ formulas, $k\ge1$, generally, to all 
existential formulas of $\sya$, as follows: 
\ben
\fenu 
\itla{modS}
If $\vpi(x)$ is a formula of $\sya$ with the only 
free variable $x$ (over $\omm$), and $T\in\pet$, 
then $T\bfo{\sus x\:\vpi(x)}$  
\imar{modS}
iff there exists $c\in\fpt$ such that $T\for\vpi(\ovc)$. 
%\qed
\een
But if $\vpi$ is a $\cip1k$ formula then we define: 
$T\bfo\vpi\,$ iff $\,T\for\vpi$.
\edf

It is a known property of the Sacks forcing that any real $x$ in 
the \dd\pet generic extension $\rV[G]$ of the universe $\rV$ has the 
form $x= \fk c(a_G)$, where $c\in\fpt\cap\rV$, see, e.g., 
\cite{nwf}. 
Therefore the forcing relation $\bfo$ as in Definition~\ref{sforc} 
is still adequate.
In particular the following lemma holds:

\ble
\lam{bfo=s}
Suppose that\/ $\vpi$ is a closed formula in\/ $\cip1k$, $k\ge1$, 
and\/ $T\in\pet$.
Then\/ $T\bfo\vpi$ iff there is no\/ $S\in\pet\yt S\sq T$, such that\/ 
$S\bfo\otr\vpi$.\qed
\ele

Here $\otr\vpi$ is the result of the canonical transformation of 
\index{zzphineg@$\otr\vpi$}
$\neg\:\vpi$ to a $\cis1k$ form.

\vyk{
\bte
\lam{bfo=s}
Suppose that\/ $\vpi$ is a closed formula in 
one of the classes\/ $\cis1k$ or\/ $\cip1k$, $k\ge1$.
Let\/ $G\sq\pet$ be Sacks-generic over the ground universe\/ $\rV$.
Then\/ $\vpi[a_G]$ is true in\/ $\rV[G]$ iff there 
is\/ $T\in G$ such that\/ $T\bfo\vpi$ in\/ $\rV$.\qed
\ete
}

Now let us address the descriptive complexity of $\bfo$. 

\ble
%[in $\rL$]
\lam{l7}
The relation\/ $\bfo$  restricted to\/ $\cip11$ 
formulas  is\/ $\hcp{1}$. 
If\/ $k\ge 2$ then the relation\/ $\bfo$ restricted to\/ $\cis1k$ 
formulas is\/ $\hcs{k-1}$ while the relation\/ $\bfo$  restricted 
to\/ $\cip1k$ formulas  is\/ $\hcp{k-1}$. 
% provided\/ $k\ge2$ and still\/ $\hcs1$ provided\/ $k=1$.
\ele
\bpf
We argue by induction.
Suppose that $\vpi=\vpi(\ovc_1,\dots,\ovc_m)$ is a closed formula 
in $\cip11$. 
It follows from the Shoenfield absoluteness and the 
perfect set theorem for $\fs11$ sets, that for any $T\in\pet$, 
$T\bfo\vpi$ is equivalent to countability of the set 
$T_\vpi=\ens{a\in[T]}{\neg\:\vpi[a]}$ in the ground universe, 
and then to 
$$
\kaz a\in[T]\:
(\vpi[a]\lor a\in\id11(c_1,\dots,c_m))
$$
as any countable set $X\sq\omm$ of class $\is11(c)$ consists
of elements of class $\id11(c)$. 
Yet the displayed formula is $\ip11$, hence $\hcd1$, 
as $x\in\id11(c)$ is a $\ip11$ relation.

The step $\ip1k\to\is1{k+1}$: make use of 
Definition~\ref{sforc}\ref{modS}. 

Now the step $\is1k\to\ip1{k}$. 
Suppose that $k\ge2$, $\vpi$ is a closed formula in $\cis1k$, and 
$T\in\pet$. 
Then by Lemma~\ref{bfo=s} $T\bfo\vpi$ is equivalent to 
$$
\kaz S\in\pet\:(S\sq T\imp \neg\;S\bfo\otr\vpi)\,,
$$ 
and hence we get $\ip1k$ using the inductive hypothesis for 
$\otr\vpi$.
\epf

\cha{Back to the \dd ncomplete Jensen's forcing}
\las{bjf} 

Let $n$ and ${\dJ}$ be as in Theorem~\ref{mt2+}.
We begin with the following

\bcl
[in $\rL$]
\lam{lcor}
For any closed formula\/ $\vpi$ in\/ $\cis1k\yt 1\le k\le{n-1}$, 
the set of all\/ 
$T\in {\dJ}$ such that\/ $T\bfo\vpi$ or\/ $T\bfo\otr\vpi$, 
is dense in\/ ${\dJ}$. 
\ecl
\bpf 
The set $W=\ens{T\in\pet}{T\bfo\vpi}$ is $\hcs{n-2}$ by 
Lemma~\ref{l7}. 
Therefore the set ${\dJ}\cap\sol W{}$ is dense in ${\dJ}$ 
by Lemma~\ref{dd}. 
However it follows from Lemma~\ref{bfo=s} that $\sol W{}$ is the 
set of all $T\in\pet$ such that $T\bfo\vpi$ or $T\bfo\otr\vpi$.
\epf

It is a basic fact with forcing that the truth in generic extensions
is certain way connected with 
the forcing relation. 
Thus the truth in \dd {\dJ}generic extensions 
$\rL[G]$ of $\rL$ corresponds to the \dd {\dJ}forcing relation. 
However, and this is the key moment, the following theorem shows
that the truth in \dd {\dJ}generic extensions is 
also in tight connection with $\pet$, the Sacks forcing notion, 
up to the level $\is1n$. 
This is a consequence of \dd ncompleteness, of course: 
in some sense, the \dd ncompleteness means that ${\dJ}$ is 
an {\it elementary submodel\/} 
of $\pet$ with respect to formulas of a certain level of 
complexity. 
%But let us go into details.

\bte
\lam{t5}
Let $n$ and ${\dJ}$ be as in Theorem~\ref{mt2+}.
Suppose that\/ $\Phi$ is a closed formula in\/ 
$\cis1k\yt 1\le k\leq n$, 
or\/ $\cip1k\yt 1\le k<n$ and 
a set\/ $G\sq {\dJ}$ is\/ \dd {\dJ}generic over\/ $\rL$.
Then $\Phi[a_G]$ holds in\/ $\rL [G]$ iff there is\/ 
$T\in G$ such that\/ $T\bfo \Phi$.
\ete
\bpf 
We argue by induction on $k$. 
Let $\Phi$ be a closed $\cip11$ formula. 
If $T\in G$ and $T\bfo\Phi$ then, in $\rL$, $\Phi[a]$ is true for 
all $a\in[T]$ with at most countable set of exceptions, see the 
proof of Lemma~\ref{l7}. 
And all exceptions are $\id11$, hence absolutely defined. 
Therefore (we skip details) 
the generic real $a_G\in[T]$ cannot be an exception, thus  
$\Phi[a_G]$ holds in $\rL [G]$. 
If $\Phi$ is $\cis11$ then $\Phi$ is $\sus x\:\vpi(x)$, $\vpi$ 
being $\cip10$, and if $T\bfo\Phi$ then by definition 
$T\bfo\vpi(\ovc)$ for some $c\in\fpt$, and so on.
On the other hand, it follows from Claim~\ref{lcor} that there 
is $T\in G$ such that $T\bfo\Phi$ or $T\bfo\otr\Phi$.
This easily implies the result for $\cis11\cup\cip11$.

\rit{Step $\cis1k\to\cip1k\yt 2\le k<n$}. 
Let $\Phi$ be a $\cip1k$ formula. 
Suppose that $\Phi[a_G]$ fails in $\rL[a_G]$. 
Then $\otr\Phi[a_G]$ holds in $\rL[a_G]$, and hence by the inductive 
hypothesis there is a condition $S\in G$ satisfying $S\bfo\otr\Phi$. 
Then by Lemma~\ref{bfo=s} there is no $T\in G$ 
with $T\bfo\Phi$. 
Conversely suppose that there is no $T\in G$ with $T\bfo\Phi$.  
Then by Claim~\ref{lcor} there is a condition 
$S\in G\yt S\bfo\otr\Phi$. 
It follows that $\otr\Phi[a_G]$ holds in $\rL[a_G]$, 
and subsequently
$\Phi[a_G]$ fails, as required. 

\rit{Step $\cip1k\to\cis1{k+1}\yt 1\le k<n$}. 
Thus let $\Phi$ be a formula $\sus x\:\vpi(x)$, where $\vpi$ is 
$\ip1k$. 
Assume that $T\in G$ satisfies $T\bfo\Phi$. 
This implies, by \ref{modS} of Definition~\ref{sforc}, that 
$T\bfo\vpi(\ovc)$ for some $c\in \fpt\cap \rL$, a code of a 
continuous map $F=\fk c:\dD\to \omm$. 
Now apply the induction hypothesis to the formula $\vpi(\ovc)$: 
it says that $\vpi(\ovc)[a_G]$ holds in $\rL[G]$. 
But $\vpi(\ovc)[a_G]$ is $\vpi[a_G](x)$, 
where $x=\fk c(a_G)\in\omm\cap\rL[G]$.
Therefore $\Phi[a_G]$ holds in $\rL[G]$, as required.

Assume, in the opposite direction, that $\Phi[a_G]$ is true in 
$\rL[G]$, that is, $\vpi[a_G](x)$ is true for some 
$x\in\omm\cap\rL[G]$. 
By Lemma~\ref{l2}, there exists 
$c\in\fpt\cap \rL$ such that $x=\fk c(a_G)$. 
The formula $\vpi(\ovc)[a_G]$ coincides with $\vpi[a_G](x)$ 
and hence 
holds in $\rL[G]$. 
Therefore, by the induction hypothesis, there is $T\in G$ 
such that 
$T\bfo \vpi(\ovc)$. 
But then $T\bfo\Phi$ by \ref{modS} of Definition~\ref{sforc}, 
as required.
\epf

\cha{Proof of the main result}
\las{pt5}

Here we accomplish the proof of Theorems \ref{mt2+} and 
\ref{mt2}.

Let $n\ge3$ and ${\dJ}\in\rL$ be as in Theorem~\ref{mt2+}.
If a set $G\sq {\dJ}$ is \dd {\dJ}generic over $\rL$ then all 
$\is1n$ sets $x\sq\om$ in $\rL[G]$ are constructible by 
Theorem~\ref{t5} because, by the homogeneity of the Sacks 
forcing, for any parameter-free formula $\Phi$ and any pair of 
trees $T\yi T'\in\pet$, it is true that 
\dm
T\bfo\Phi\;\eqv\; S\bfo\Phi. 
\dm
Let us present this final argument in more detail. 

If $S,T\in\pet$ then let $\Hom_{ST}$ be the set of all 
homeomorphisms $h:[S]\onto[T]$; clearly $\Hom_{ST}$ is non-empty. 
Suppose that $h\in\Hom_{ST}$. 
Recall that continuous functions $F:\dD\to\omm$ are coded so that 
$\fk c$ is the function coded by $c\in\fpt$. 
If $c,d\in\fpt$ then write $c\equ h d$ iff 
$\fk d(h(a))= \fk c(a)$ for all $a\in[S]$.\snos
{We suppress the dependence on $S$ in $c\equ h d$ 
since formally $S=\dom h$.}
If $\vpi=\Phi(\ovc_1,\dots,\ovc_m)$ and $\psi=\Phi(\ovd_1,\dots,\ovd_m)$ 
are formulas of $\sya$ (see \nrf{dsf}), 
and $c_i\equ h d_i$ for all $i$, then write $\vpi\equ h\psi$. 
In this case the formulas 
$\vpi[a]$ and $\psi[h(a)]$ coincide for any $a\in[S]$.

\bcl
\lam{he2}
Suppose that\/ $S,T\in\pet$, $h\in\Hom_{ST}$, $\Phi$ and\/ 
$\Psi$ are closed formulas in one and the same class\/ $\cis1k$ 
or\/ $\cip1k$, and\/ $\Phi\equ h\Psi$.
Then\/ $S\bfo\Phi$ if and only if\/ $T\bfo\Psi$.
\ecl
\bpf
Exercise, argue by induction on the complexity of the formulas.
\epf

\bcor
\lam{he3}
If\/ $S,T\in\pet$  
and\/ $\Phi$ is a closed formula in\/ $\is1k$ or\/ $\ip1k$ 
%not containing any symbols\/ $\ovc$, 
then\/ $S\bfo\Phi$ if and only if\/ $T\bfo\Phi$.
\ecor
\bpf
Pick $h\in\Hom_{ST}$, note that $\Phi\equ h\Phi$ 
(as $\Phi$ contains no symbols of the form $\ovc$), 
and apply Claim~\ref{he2}. 
\epf

Now the next lemma completes the proof of Theorem~\ref{mt2}.

\ble
\lam{he4}
If\/ $G\sq\rL[G]$ is\/ \dd {\dJ}generic over\/ $\rL$ and\/ 
$x\sq\om\yt x\in\rL[G]$ is\/ $\is1n$ in\/ $\rL[G]$ then\/  
$x\in\rL$ and\/ $x$ is\/ $\is1n$ in\/ $\rL$. 
\ele
\bpf
Let $\vpi(m)$ be a parameter-free $\is1n$ formula such that 
$x=\ens{m}{\vpi(m)}$ in $\rL[G]$.
Consider the tree $S=\bse\in\pet$. 
Then
$$
m\in x\leqv \rL[G]\models\vpi(m)\leqv
\sus T\in G\:(T\bfo\vpi(m))\leqv
S\bfo\vpi(m)\,,
$$
where the two last equivalences follow from resp.\ Theorem~\ref{t5} 
and Corollary~\ref{he3}.
It remains to apply Lemma~\ref{l7}.
\epf

\qeD{Theorems \ref{mt2+} and \ref{mt2}}

\cha{Further applications}
\las{fapp}

Countable support products of \rit{independent} 
forcing notions of Jensen's type have been recently 
applied for some definability problems. 
Let us review some of the applications. 

1. 
A somewhat modified forcing notion, 
say ${\dJ}'\sq\pet$, 
rather similar to the above forcing ${\dJ}$, 
is defined in \cite{kl27e}.  
Its finite support product $(\dJ')^\om$ 
still satisfies CCC in $\rL$,  
and if $G\sq(\dJ')^\om$ is generic 
over $\rL$ then it is true in $\rL[G]$ that 
the corresponding set 
$\ens{x_k[G]}{k<\om}$ 
of \dd{\dJ'}generic reals $x_k[G]$ 
is a countable lightface $\ip12$ set containing 
no ordinal-definable elements. 

2. 
Another, more profound modification $\dJ^{\text{inf}}$ 
of the above forcing ${\dJ}$ 
is defined in \cite{kl22}.  
It still satisfies CCC in $\rL$, 
and is invariant under some simple transformations, 
so, instead of a single generic real as $\dJ$ does, 
it adds a \dd\Eo\rit{equivalence class}\snos
{Two reals $a,b\in\dn$ are \dd\Eo equivalent if 
$a(n)=b(n)$ for all but finite $n$.} 
of \dd{\dJ^{\text{inf}}}generic reals, 
and it turns out that this \dd{\dJ^{\text{inf}}}generic 
\dd\Eo class is a (countable) $\ip12$ set containing 
no OD (ordinal-definable) elements in the extension.
A related forcing notion in \cite{kl25} yields a 
generic extension of $\rL$ in which there exist 
two different OD-indiscernible \dd\Eo equivalence 
classes, whose union is lightface $\ip12$. 
See \cite{fgh} on modern development around 
OD-indiscernibility and related notions. 

3. 
A generalization of $\dJ^{\text{inf}}$, 
obtained in \cite{kl34}, provides a model, in which, 
for a given $n\ge2$, there is a countable $\ip1n$ 
set containing no ordinal-definable elements and 
equal to a \dd\Eo class, and in the same time 
every countable $\is1n$ 
set definitely contains  ordinal-definable elements.
 
4. 
A sequence\/ $\sis{{\dJ}_\xi}{\xi<\omi}\in\rL$ of forcing 
notions\/ ${\dJ}_\xi\sq\pet$, rather similar to the 
above forcing ${\dJ}$, is defined in \cite{kl28}.  
The finite support product $\dJ^{\omi}=\prod_\xi\dJ_\xi$ 
satisfies CCC in $\rL$,  
``being a \dd{{\dJ}_\xi}generic real'' 
is\/ $\ip\hc{1}$ uniformly in\/ $k$, 
and in addition the forcing notions\/ ${\dJ}_\xi$ 
are independent in the 
sense that for any $\xi<\omi$ there is no 
\dd{{\dJ}_\xi}generic reals 
in \dd{\prod_{\et\ne\xi}{\dJ}_\eta}generic extensions 
of $\rL$. 
Different subextensions of \dd{\dJ^{\omi}}generic 
extensions of $\rL$ have interesting definability 
effects, for instance, the failure of the separation  
theorem for $\fp13$, or the existence of a 
planar $\ip12$ set with countable sections, 
which cannot be uniformized by a ROD 
(real-ordinal definable) 
set, \cite{kl28}.  

5.
The last result was generalized in \cite{kl38}: 
a model, in which, 
for a given $n\ge2$, there is a planar 
ROD-nonuniformizable $\ip1n$ set with countable 
cross-sections, and in the same time 
every planar $\fs1n$ set with countable 
cross-sections is 
ROD-uniformizable.

6.
Even the countable subsequence 
$\sis{{\dJ}_k}{k<\om}\in\rL$ of the \dd\omi sequence 
of Item 4 yields an interesting result in a sense  
complementary to Theorem~\ref{mt2}.
Namely let $\vec a=\sis{a_k}{k<\om}$ be a finite-support 
\dd{\prod_{k<\om}\dJ_k}generic sequence 
over $\rL$, of reals $a_k\in\dn.$ 
Let $x\in\dn$ be Cohen-generic over $\rL[\vec a]$.
Then it is true in the submodel 
$\rL[x,\sis{a_k}{x(k)=1}]\sq \rL[\vec a]$
that $x$ is a $\is13$ set while all $\id13$ sets 
$y\sq\om$ are constructible.

7. 
Abraham \cite{abr2} developed a Jensen-like forcing 
by \dd\omi branching trees, which adjoins a generic 
collapsing map $f:\om\onto\omil$ to $\rL$, 
effectively coded by a $\id13$ real. 
An \dd nlevel generalization was obtained in 
\cite{kl36}: if $n\ge3$ then there is a generic 
extension of $\rL$ by a collapsing map 
$f:\om\onto\omil$ to $\rL$, 
effectively coded by a $\id1{n+1}$ real, 
and in the same time it is true in $\rL$ that 
no $\id1n$ real codes a collapsing map.

8. 
Iterations of Jensen's forcing were developed by 
Abraham~\cite{abr}.   
Combining this technique with finite-support 
Jensen products technique of \cite{kl22} 
and some earlier forcing constructions used in the 
theory of generic choiceless models, a model 
of $\ZF$ is constructed in \cite{kfg} in which 
the countable $\AC$ holds but the Dependent Choices 
scheme $\DC$ fails for some $\ip12$ formula 
(which is the best possible). 

9. 
We finish with a review of several marvellous 
results of Harrington, known only from his 
handwritten notes \cite{h74}, never published 
neither in full nor in part. 

\cite[Part A]{h74} 
Harrington sketches a generic extension of $\rL$ in which 
the set $\dn\cap\rL$ of all constructible reals 
is equal to the set $\id13\cap\dn$ of all 
lightface $\id13$ reals, and adds a few 
sentences to explain what should be done in order 
to get a model for $\dn\cap\rL=\id1{n+1}\cap\dn$, 
given $n\ge3$, and a model for 
$\dn\cap\rL=\id1{\infty}\cap\dn$. 
The construction involved the almost-disjoint 
coding of Jensen -- Solovay \cite{jsad}. 
The result was also mentioned in  
\cite{hf102} as a solution of Problem 87 there. 
An independent and somewhat different construction 
of a model for 
$\dn\cap\rL=\id1{\infty}\cap\dn$ 
(constructible reals = analytically definable reals), 
based on somewhat different technique, 
was given in in \cite{ian79}.
Needless to say that $\omil$ has to be countable in all 
these models.

\cite[Part B]{h74} 
A model is sketched, in which the separation theorem 
fails for both classes $\fs13$ and $\fp13$. 
The failure of \dd{\fp13}separation is especially 
interesting since the separation theorem 
holds for $\fp13$ in $\rL$. 
A hint is given, as how to construct a model in 
which the separation theorem 
fails for both classes $\fs1n$ and $\fp1n$, 
given $n\ge3$, or even immediately for all 
those classes. 
(Another model, where \dd{\fp13}separation 
fails, but based on Jensen forcing 
products rather than almost-disjoint forcing 
technique as in \cite{h74}, 
is presented in \cite{kl28}.) 
 
\cite[Part C]{h74} 
An extremely complicated construction of a model 
is sketched, in which \dd{\fp13}separation holds, 
but \dd{\fs13}reduction fails. 
Beside the almost-disjoint coding, 
the construction involves a careful study of 
subextensions of the almost-disjoint generic 
extensions considered. 
Generalization to $n\ge4$ was left as a matter of 
belief.

\cite[Part D]{h74} 
A model is sketched, in which 
the separation theorem 
holds for both lightface classes $\is13$ and $\ip13$, 
but only for sets of integers.
It is asserted that, with some changes, the result 
can be obtained for $\is1n$ and $\ip1n$, 
any given $n\ge4$.\vtm

{\bf Acknowledgements.} 
We thank the anonymous reviewer for 
their thorough review and highly appreciate the comments 
and suggestions, which significantly contributed 
to improving the quality of the publication. 

}

\renek{\refname}     {\large\bf References}

\bibliographystyle{plain}

%\vyk
{\small
%
%\bibliography{51x,51klx}
%

}

\mtho

\end{document}